\documentclass[11pt]{article}

\usepackage{amssymb}
\usepackage{amsmath}
\usepackage{amsthm}
\usepackage{mathrsfs}
\usepackage{graphicx}
\usepackage[utf8]{inputenc}
\usepackage{multirow}
\usepackage{floatrow}
\usepackage{mathrsfs}
\usepackage{xcolor}
\usepackage{import}
\usepackage{xifthen}
\usepackage{pdfpages}
\usepackage{transparent}

\font\smallit=cmti10

\newtheorem{thm}{Theorem}
\newtheorem{prop}[thm]{Proposition}

\newtheorem{cor}[thm]{Corollary}
\newtheorem{rem}[thm]{Remark}
\theoremstyle{definition}
\newtheorem{eg}[thm]{Example}

\newenvironment{Pf}{\noindent{\bf Proof. }}{\hfill $\blacksquare$ \\}

\def\ts{\textstyle}
\def\lf{\lfloor}
\def\rf{\rfloor}
\def\lc{\lceil}
\def\rc{\rceil}
\def\imp{\rightarrow}

\def\LIff{\Longleftrightarrow}
\def\la{\langle}
\def\ra{\rangle}
\def\pr{\prime}
\def\a{\alpha}
\def\b{\beta}

\def\N{\mathbb N}
\def\Q{\mathbb Q}
\def\Z{\mathbb Z}
\def\l{\lambda}
\def\sm{\setminus}

\newsavebox{\f}
\savebox{\f}{$\quad\;\left\{\text{\smallit maz188235}, \text{\smallit atripath} \right\}$}
\begin{document}

\title{\bf On proportionally modular numerical semigroups that are generated by arithmetic progressions}
\author
{{\bf Edgar Federico Elizeche}\thanks{{\smallit Department of Mathematics, Indian Institute of Technology, Hauz Khas, New Delhi -- 110016, India} \newline {\usebox{\f}\smallit @maths.iitd.ac.in}} \and
{\bf Amitabha Tripathi} \footnotemark[1]\:\:\thanks{Corresponding author}}
\date{}
\maketitle

\begin{abstract}
\noindent A numerical semigroup is a submonoid of ${\Z}_{\ge 0}$ whose complement in ${\Z}_{\ge 0}$ is finite. For any set of positive integers $a,b,c$, the numerical semigroup $S(a,b,c)$ formed by the set of solutions of the inequality $ax \bmod{b} \le cx$ is said to be proportionally modular. For any interval $[\a,\b]$, $S\big([\a,\b]\big)$ is the submonoid of ${\Z}_{\ge 0}$ obtained by intersecting the submonoid of ${\Q}_{\ge 0}$ generated by $[\a,\b]$ with ${\Z}_{\ge 0}$. 

\noindent For the numerical semigroup $S$ generated by a given arithmetic progression, we characterize $a,b,c$ and $\a,\b$ such that both $S(a,b,c)$ and $S\big([\a,\b]\big)$ equal $S$.    
\end{abstract}


\vskip 20pt

\section{Introduction} \label{intro}
\vskip 10pt

A numerical semigroup $S$ is a submonoid of ${\Z}_{\ge 0}$ whose complement ${\Z}_{\ge 0} \sm S$ is finite. For the complement to be finite, it is necessary and sufficient that $\gcd(S)=1$. For a given subset $A$ of positive integers, we write
\[ \la A \ra = \big\{ a_1x_1+\cdots+a_kx_k: a_i \in A, x_i \in {\Z}_{\ge 0}, k \in \N \big\} \:\bigcup\: \{0\}. \]
Note that $\la A \ra$ is a submonoid of ${\Z}_{\ge 0}$, and that $S=\la A \ra$ is a numerical semigroup if and only if $\gcd(A)=1$. 

For a numerical semigroup $S$, the complement ${\Z}_{\ge 0} \sm S$ is denoted by $\text{G}(S)$, and called the gap set of $S$. The largest element in $\text{G}(S)$ is called the Frobenius number of $S$, and denoted by $\text{F}(S)$. The cardinality of $\text{G}(S)$ is called the genus of $S$, and denoted by $g(S)$.   

We say that $A$ is a set of generators of the semigroup $S$, or that the semigroup $S$ is generated by the set $A$, when $S=\la A \ra$. Further, $A$ is a minimal set of generators for $S$ if $A$ is a set of generators of $S$ and no proper subset of $A$ generates $S$. Every semigroup has a unique minimal set of generators. The embedding dimension $e(S)$ of $S$ is the size of the minimal set of generators. 

Given a subset $A$ of ${\Q}_{\ge 0}$, we denote by $\la A \ra$ the submonoid of ${\Q}_{\ge 0}$ generated by $A$: 
\[ \la A \ra = \big\{ a_1x_1+\cdots+a_kx_k: a_i \in A, x_i \in {\Z}_{\ge 0}, k \in \N \big\} \:\bigcup\: \{0\}. \]
Then $S(A)=\la A \ra \cap {\Z}_{\ge 0}$ is a submonoid of ${\Z}_{\ge 0}$. We say that $S(A)$ is the numerical semigroup associated to $A$. When $A=I$ is an interval, there is a simple way to describe $\la A \ra$.   
\vskip 5pt

\begin{prop} {\bf (\cite[Lemma 1, pp. 282]{RGGU03})} \label{basic1}
Let $x \in {\Q}^+$ and let $I$ be an interval. Then $x \in \la I \ra$ if and only if there exists a positive integer $n$ such that $\frac{x}{n} \in I$. 
\end{prop}
\vskip 5pt

\noindent Toms \cite{Tom03} in an attempt to solve problems in classification theory in $C^{\star}$-algebras defined Toms decomposable numerical semigroups. The concept of proportionally modular numerical semigroups was introduced by Rosales {\it et al\/} in \cite{RGU05}. Rosales \& Garcia-Sanchez \cite{RG-S08} showed that a numerical semigroup is Toms decomposable if and only if it is the intersection of finitely many proportionally modular numerical semigroups. 
\vskip 5pt

\noindent A proportionally modular Diophantine inequality is an expression of the form $ax \bmod{b} \le cx$, where $a,b,c$ are positive integers. The set of integer solutions of a proportionally modular Diophantine inequality form a numerical semigroup. For $a,b,c \in \N$, define
\[ S(a,b,c) = \{ x \in {\Z}_{\ge 0}: ax \bmod{b} \le cx \}. \]
A numerical semigroup of this form is called proportionally modular. Since the inequality $ax \bmod{b} \le cx$ has the same set of integer solutions as $(a \bmod{b})x \bmod{b} \le cx$, we may assume that $a \in \{1,\ldots,b\}$. We note that $S(b,b,c)={\Z}_{\ge 0}$. 
\vskip 5pt

\noindent Given a proportionally modular numerical semigroup, there exists a dual proportionally modular numerical semigroup. 

\begin{prop} {\bf (\cite[Proposition 1, pp. 416]{DR06})} \label{basic2}
Let $a,b,c \in \N$ with $c<a<b$. Then 
\[ S(a,b,c) = S(b+c-a,b,c). \] 
\end{prop}
\vskip 5pt

We note that $S\big([\a,\b]\big)={\Z}_{\ge 0}$ if and only if $\a \le \frac{1}{n} \le \b$ for some $n \in \N$, and that $S(a,b,c)={\Z}_{\ge 0}$ if and only if $a \bmod{b} \le c$. The following two Propositions give a connection between $S(a,b,c)$ and $S\big([\a,\b]\big)$. 
\vskip 5pt 

\begin{prop} {\bf (\cite[Lemma 1, pp. 454]{RGU08})} \label{basic3}
Let $a,b,c \in \N$ with $c<a$. Then 
\[ S(a,b,c) = S\left(\left[\ts\frac{b}{a},\ts\frac{b}{a-c}\right]\right). \] 
\end{prop}

\begin{prop} {\bf (\cite[Lemma 1, pp. 454]{RGU08})} \label{basic4} 
Let $\frac{p}{q},\frac{r}{s} \in \Q$ with $\frac{p}{q}<\frac{r}{s}$. Then 
\[ S\left(\left[\ts\frac{p}{q},\ts\frac{r}{s}\right]\right) = S(qr,pr,qr-ps). \] 
\end{prop}
\vskip 5pt

Let $a,b,c \in \N$ with $c<a<b$. Using Propositions \ref{basic2} and \ref{basic3}, there exist intervals $[\a,\b]$ and $[{\a}^{\pr},{\b}^{\pr}]$, with rational endpoints, and with $\a>1$ and ${\a}^{\pr}>1$, such that 
\[ S\left([\a,\b]\right) = S\left([{\a}^{\pr},{\b}^{\pr}]\right). \] 
Moreover, it easily follows from Propositions \ref{basic2} and \ref{basic3} that 
\begin{equation} \label{basic5}
{\a}^{\pr} = \frac{\b}{\b-1}, \quad {\b}^{\pr} = \frac{\a}{\a-1}. 
\end{equation}
\vskip 5pt

\noindent B\'{e}zout sequences are closely connected to the study of proportionately modular numerical semigroups; see \cite{RGU08}. A B\'{e}zout sequence is an increasing finite sequence of rational numbers $\left\{\frac{p_1}{q_1},\ldots,\frac{p_n}{q_n}\right\}$ with $\gcd(p_k,q_k)=1$ satisfying $p_{k+1} q_k-p_k q_{k+1}=1$ for $k \in \{1,\ldots,n-1\}$.  A B\'{e}zout sequence is said to be {\em proper} if $p_{k+\ell} q_k-p_{k+\ell} q_k>1$ for each $k$ and each $\ell>1$. As a consequence, we have 
\begin{equation} \label{convex} 
p_1 > \cdots > p_r < p_{r+1} < \cdots < p_n 
\end{equation}
for some $r \in \{1,\ldots,n\}$; see \cite[Corollary 18, pp. 459]{RGU08}.    

From the proof of Lemma 22 and Theorem 23 in \cite{RGU08}, given an interval $I=[\a,\b]$, $\a,\b \in \Q$ and $S=S(I)$, there exists a permutation of the generators $a_1,\ldots,a_e$ of $S$ and positive integers $b_1,\ldots,b_e$ such that 
\begin{equation} \label{Bez1}
\left\{\frac{a_n}{b_n}\right\}_{n \ge 1} \:\text{forms a proper B\'{e}zout sequence, and}\:\a \le \frac{a_1}{b_1} < \frac{a_e}{b_e} \le \b. 
\end{equation}  
\vskip 5pt

\noindent Given rational numbers $p_1/q_1$ and $p_n/q_n$ in reduced form, there exists a unique proper B\'{e}zout sequence with first term $p_1/q_1$ and last term $p_n/q_n$. This unique B\'{e}zout sequence has been determined algorithmically; see \cite[Algorithm 3.5]{BR09}.   
\vskip 5pt

The following Proposition gives a characterization of proportionally modular numerical semigroups in terms of its minimal set of generators. 
\vskip 5pt

\begin{prop} {\bf (\cite[Theorem 31, pp. 463]{RGU08})} \label{Bez2} 
A numerical semigroup $S$ is proportionally modular if and only if there exists a permutation $a_1,\ldots,a_e$ of its minimal set of generators such that the following two conditions hold: 
\begin{itemize}
\item[{\rm (i)}]  
$\gcd(a_k,a_{k+1})=1$ for $1 \le k \le e-1$, and \\[-16pt]
\item[{\rm (ii)}]
$a_k \mid (a_{k-1}+a_{k+1})$ for $2 \le k \le e-1$. 
\end{itemize}  
\end{prop}
\vskip 5pt

\noindent Certain aspects of the numerical semigroup generated by an arithmetic progression $a,a+d,a+2d,\ldots,a+(k-1)d$ has been studied by many authors; see \cite{Bat58, Bra42, Gra73, NW72, Rob56, Tri94}. In particular, the gap set $\text{G}(S)$ and the Frobenius number $\text{F}(S)$ are determined by the following result in \cite{Tri94}. It can be shown that the embedding dimension $e(S)=\min\{a,k\}$. 

\begin{prop} {\bf (\cite[Theorem 1, pp. 780]{Tri94})} \label{T} 
Let $a,d,k$ be positive integers, with $\gcd(a,d)=1$ and $k \ge 2$. Let $S=\la a,a+d,a+2d,\ldots,a+(k-1)d \ra$. Then 
\begin{itemize}
\item[{\rm (i)}]  
$\text{G}(S) = \left\{ ax+dy: 0 \le x \le \left\lf \ts\frac{y-1}{k-1} \right\rf, 1 \le y \le a-1 \right\}$. \\[-12pt]
\item[{\rm (ii)}]
$\text{F}(S) = \max \text{G}(S) = a\left\lf \ts\frac{a-2}{k-1} \right\rf+d(a-1)$. 
\end{itemize}  
\end{prop}
\vskip 5pt

\noindent Given a numerical semigroup $S$, it is natural to ask for a characterization of positive integers $a,b,c$ for which $S(a,b,c)=S$, and also for a characterization of intervals $I=[\a,\b]$ for which $S\big([\a,\b]\big)=S$. We note that Propositions \ref{basic2} and \ref{basic3} connect the numerical semigroups $S(a,b,c)$ and $S\big([\a,\b]\big)$, so characterizing one characterizes the other. Let $S=S\big(AP(a,d;k)\big)=\la a,a+d,a+2d,\ldots,a+(k-1)d \ra$, $k \ge 2$. Since $S\big(AP(a,d;k)\big)=S\big(AP(a,d;a)\big)$ whenever $k>a$, we may assume, without loss of generality, that $k \le a$. 

In this paper, we characterize $a,b,c$ and $\a,\b$ such that both $S(a,b,c)$ and $S\big([\a,\b]\big)$ equal the numerical semigroup $S=S\big(AP(a,d;k)\big)$ generated by any arithmetic progression. The two characterizations appear as Theorems \ref{I} and \ref{II}, and are based on Propositions \ref{max}, \ref{min}, and \ref{link}. Delgado \& Rosales \cite{DR06} provide an algorithmic formula to determine $\text{F}\big(S(a,b,c)\big)$ that does not lead to an explicit formula. The determination of an explicit formula for $\text{F}\big(S(a,b,c)\big)$ therefore remains an open problem. Our characterization of $a,b,c$ for which $S=S\big(AP(a,d;k)\big)$ yields an explicit formula in these special cases since the problem of determining the Frobenius number for semigroups generated by arithmetic progressions is well known, as mentioned above. We note in passing that the algorithm in \cite[Algorithm 3.5]{BR09} to compute the unique proper B\'{e}zout sequence connecting two given reduced rational numbers $p_1/q_1$ and $p_n/q_n$ relies on knowing these rational numbers. Since there is no obvious relationship between $p_1/q_1$ and $p_n/q_n$, and the numerical semigroup $S$, these endpoints of the proper B\'{e}zout sequence has to be determined first. This means that the algorithm in \cite[Algorithm 3.5]{BR09} is not directly useful in resolving the problem.   
 
\section{Main Results} \label{main}
\vskip 10pt

\noindent In this section, we consider an inverse problem related to two types of numerical semigroups that were introduced in Section \ref{intro} -- numerical semigroups generated by intervals $[\a,\b]$ and proportionally modular numerical semigroups $S(a,b,c)$. Propositions \ref{basic2} and \ref{basic3} in Section \ref{intro} connect these two types of numerical semigroups. 

In order that the semigroup generated by an arithmetic progression $a,a+d,a+2d,\ldots,a+(k-1)d$ be a numerical semigroup, it is necessary and sufficient that $\gcd(a,d)=1$. Henceforth we assume $\gcd(a,d)=1$, and write $S\big(AP(a,d;k)\big)=\la a,a+d,a+2d,\ldots,a+(k-1)d \ra$. The set $\{a,a+d,a+2d,\ldots,a+(k-1)d\}$ is the minimal set of generators for $S\big(AP(a,d;k)\big)$ when $2 \le k \le a$, whereas the set $\{a,a+d,a+2d,\ldots,a+(a-1)d\}$ is the minimal set of generators for $S\big(AP(a,d;k)\big)$ for $k>a$. Throughout this section, we may therefore assume without loss of generality that $2 \le k \le a$. 
\vskip 5pt

\noindent We illustrate our main result, Theorem \ref{I}, with the following two examples. In Example \ref{SG1} we determine intervals $I$ for which $S(I)$ is the ordinary numerical semigroup $\{0,5,\imp\}$, and in Example \ref{SG2} intervals $I$ for which $S(I)$ is the numerical semigroup $\la 5,7 \ra$ with embedding dimention two. These are both special cases of our result in Theorem \ref{I}.  
\vskip 5pt

\begin{eg} \label{SG1}
We find $I=[\a,\b]$, $\a>1$, for which $S(I)=\{0,5,\imp\}$. From the argument leading up to eqn.~\eqref{Bez1}, we have $\{a_1,\ldots,a_5\}=\{5,\ldots,9\}$, giving rise to the following four B\'{e}zout sequences: 
\begin{itemize}
\item $\frac{5}{1}$, $\frac{6}{1}$, $\frac{7}{1}$, $\frac{8}{1}$, $\frac{9}{1}$. 
\item $\frac{9}{2}$, $\frac{5}{1}$, $\frac{6}{1}$, $\frac{7}{1}$, $\frac{8}{1}$.  
\item $\frac{9}{8}$, $\frac{8}{7}$, $\frac{7}{6}$, $\frac{6}{5}$, $\frac{5}{4}$. 
\item $\frac{8}{7}$, $\frac{7}{6}$, $\frac{6}{5}$, $\frac{5}{4}$, $\frac{9}{7}$.  
\end{itemize}
This leads to  
\[ \a \le \tfrac{5}{1} < \tfrac{9}{1} \le \b \;\;\text{or}\;\; \a \le \tfrac{9}{8} < \tfrac{5}{4} \le \b \;\;\text{or}\;\; \a \le \tfrac{9}{2} < \tfrac{8}{1} \le \b \;\;\text{or}\;\; \a \le \tfrac{8}{7} < \tfrac{9}{7} \le \b. \] 
Note that the pairs of intervals $[5,9]$, $[\tfrac{9}{8},\tfrac{5}{4}]$ and $[\tfrac{9}{2},8]$, $[\tfrac{8}{7},\tfrac{9}{7}]$ are connected by eqn.~\eqref{basic5}. 
\vskip 5pt

\noindent From Proposition \ref{basic1}, $\frac{g}{n} \notin I$ whenever $g \in \text{G}(S)=\{1,2,3,4\}$ and $n \in \N$. This results in 
\[ \a \in (4,5], \b \in [9,\infty) \;\;\text{or}\;\; \a \in \left(1,\tfrac{9}{8}\right], \b \in \left[\tfrac{5}{4},\tfrac{4}{3}\right) \;\;\text{or}\;\; 
    \a \in \left(4,\tfrac{9}{2}\right], \b \in [8,\infty) \;\;\text{or}\;\; \a \in \left(1,\tfrac{8}{7}\right], \b \in \left[\tfrac{9}{7},\tfrac{4}{3}\right).
\]
Note that the pairs of intervals $(4,5]$, $\left[\tfrac{5}{4},\tfrac{4}{3}\right)$, $[9,\infty)$, $\left(1,\tfrac{9}{8}\right]$,  $\left(4,\tfrac{9}{2}\right]$, $\left[\tfrac{9}{7},\tfrac{4}{3}\right)$, and $[8,\infty)$, $\left(1,\tfrac{8}{7}\right]$ are connected by eqn.~\eqref{basic5}. 

\begin{figure}[h!] 
\def\svgwidth{\columnwidth}
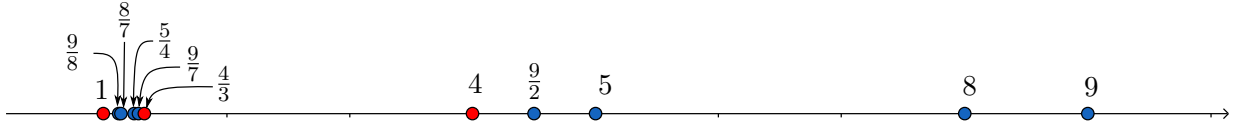

\caption{For $I$ such that $S(I)=\{0,5,\imp\}$, $[5,9] \subset I \subset (4,\infty)$ or $[\frac{9}{2},8] \subset I \subset (4,\infty)$ or $[\frac{9}{7},\frac{5}{4}] \subset I \subset (1,\frac{4}{3})$ or $[\frac{8}{7},\frac{9}{7}] \subset I \subset (1,\frac{4}{3})$.}
\end{figure}
\end{eg}
\vskip 5pt

\begin{eg} \label{SG2}
From Proposition \ref{Bez2}, it follows that every numerical semigroup with embedding dimension $2$ is proportionally modular. We find $I=[\a,\b]$, $\a>1$, for which $S(I)=\la 5,7 \ra$. From the argument leading up to eqn.~\eqref{Bez1}, we have $\{a_1,a_2\}=\{5,7\}$, so that $1<\frac{5}{b_1}<\frac{7}{b_2}$ and $7b_1-5b_2=1$ or $1<\frac{7}{b_1}<\frac{5}{b_2}$ and $5b_1-7b_2=1$. This leads to  
\[ \a \le \tfrac{5}{3} < \tfrac{7}{4} \le \b \;\;\text{or}\;\; \a \le \tfrac{7}{3} < \tfrac{5}{2} \le \b. \] 
Note that the intervals $[\tfrac{5}{3},\frac{7}{4}]$ and $[\tfrac{7}{3},\tfrac{5}{2}]$ are connected by eqn.~\eqref{basic5}. 
\vskip 5pt

\noindent From Proposition \ref{basic1}, $\frac{g}{n} \notin I$ whenever $g \in \text{G}(S)=\{1,2,3,4,6,8,9,11,13,16,18,23\}$ and $n \in \N$. This results in 
\[ \a \in \left(\tfrac{23}{14},\tfrac{5}{3}\right], \b \in \left[\tfrac{7}{4},\tfrac{23}{13}\right) \quad \text{or} \quad 
    \a \in \left(\tfrac{23}{10},\tfrac{7}{3}\right], \b \in \left[\tfrac{5}{2},\tfrac{23}{9}\right). 
\]
Note that the pairs of intervals $\left(\tfrac{23}{14},\tfrac{5}{3}\right]$, $\left[\tfrac{5}{2},\tfrac{23}{9}\right)$ and $\left[\tfrac{7}{4},\tfrac{23}{13}\right)$, $\left(\tfrac{23}{10},\tfrac{7}{3}\right]$ are connected by eqn.~\eqref{basic5}. 

\begin{figure}[h!] 
\def\svgwidth{\columnwidth}
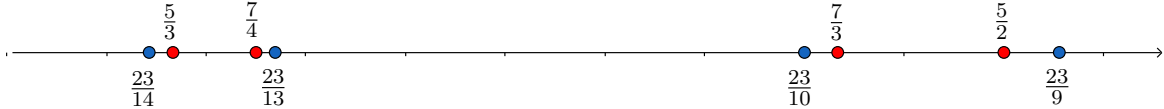

\caption{For $I$ such that $S(I)=\la 5,7 \ra$, $[\frac{5}{3},\frac{7}{4}] \subset I \subset (\frac{23}{14},\frac{23}{13})$ or $[\frac{7}{3},\frac{5}{2}] \subset I \subset (\frac{23}{10},\frac{23}{9})$.}
\end{figure}
\end{eg}
\vskip 5pt

The main result in this section is the characterization of $\a,\b$ such that $S\big([\a,\b]\big)=S\big(AP(a,d;k)\big)$ in Theorem \ref{I}. Essential to our proof of Theorem \ref{I} are three results listed as Propositions \ref{max}, \ref{min}, and \ref{link}. We also characterize $a,b,c$ such that $S(a,b,c)=S\big(AP(a,d;k)\big)$, via Proposition \ref{basic2} and Theorem \ref{I}. We also state, without proof, the two special cases $S\big(AP(a,1;k)\big)$ and $S\big(AP(a,1;a)\big)$.  
\vskip 5pt

\begin{prop} \label{max}
Let $a,d,k$ be positive integers, with $\gcd(a,d)=1$ and $2 \le k \le a$. Let $S\big(AP(a,d;k)\big)=\la a,a+d,a+2d,\ldots,a+(k-1)d \ra$. If $dd_a=1+q_a a$, with $1 \le d_a<a$ and $0 \le q_a<d$, and $\l=\left\lf\frac{a-2}{k-1}\right\rf$, then 
\[ \max \left\{ \frac{g}{\left\lc \frac{d_a g}{a} \right\rc}: g \in \text{G}(S) \right\} =  \frac{\l a+d(a-1)}{\l d_a+q_a(a-1)+1} = 
     \frac{\text{F}(S)}{\left\lc \frac{d_a\text{F}(S)}{a} \right\rc}, 
\]
with the maximum achieved at $g=\l a+d(a-1)$. 
\end{prop}

\begin{Pf}
If $g \in \text{G}(S)$, then $g=ax+dy$, with $0 \le x \le \left\lf \frac{y-1}{k-1} \right\rf$ and $1 \le y \le a-1$. Let $d_a$ be such that $dd_a \equiv 1\pmod{a}$, $d_a>0$, and write $dd_a=q_a a+1$, $q_a \in {\Z}_{\ge 0}$. Then 
\[ \frac{d_a g}{a} = \frac{d_a(ax+dy)}{a} = d_a x + q_a y + \frac{y}{a}. \]
Since $1 \le y<a$, 
\[ \frac{g}{\left\lc \ts\frac{d_a g}{a} \right\rc} = \frac{ax+dy}{d_a x+q_a y+1} \stackrel{\mbox{\tiny{def}}}{=} f_{\ell}(x,y). \]
We show that 
\begin{equation} \label{lbasic1}
\max \left\{ f_{\ell}(x,y): 0 \le x \le \left\lf \ts\frac{y-1}{k-1} \right\rf, 1 \le y \le a-1 \right\} = f_{\ell}\left( \left\lf \ts\frac{a-2}{k-1} \right\rf, a-1\right) = \frac{\text{F}(S)}{\left\lc \frac{d_a\text{F}(S)}{a} \right\rc}. 
\end{equation}
\vskip 5pt

\noindent Fix $y \in \{1,\ldots,a-1\}$. For $x \in \left\{0,1,2,\ldots,\left\lf \ts\frac{y-1}{k-1} \right\rf -1\right\}$, the numerator of 
\[ f_{\ell}(x+1,y) - f_{\ell}(x,y) = \frac{a(x+1)+dy}{d_a(x+1)+q_a y+1} - \frac{ax+dy}{d_a x+q_a y+1} \]
is 
\[ a(d_a x+q_a y+1)-d_a(ax+dy) = a(q_a y+1)-y(q_a a+1) =  a-y > 0. \]
Therefore, 
\begin{equation} \label{lbasic2}
\max \left\{ f_{\ell}(x,y): 0 \le x \le \left\lf \ts\frac{y-1}{k-1} \right\rf, 1 \le y \le a-1 \right\} = \max \left\{ f_{\ell}\left(\left\lf \ts\frac{y-1}{k-1} \right\rf,y\right): 1 \le y \le a-1 \right\}. 
\end{equation}
\vskip 5pt

\noindent For $y \in \{1,\ldots,a-2\}$, the numerator of 
\[ f_{\ell}\left(\left\lf \ts\frac{y}{k-1} \right\rf,y+1\right) - f_{\ell}\left(\left\lf \ts\frac{y-1}{k-1} \right\rf,y\right) =  
    \frac{a\left\lf \ts\frac{y}{k-1} \right\rf+d(y+1)}{d_a\left\lf \ts\frac{y}{k-1} \right\rf+q_a(y+1)+1} - 
    \frac{a\left\lf \ts\frac{y-1}{k-1} \right\rf+dy}{d_a\left\lf \ts\frac{y-1}{k-1} \right\rf+q_a y+1} 
\]
is 
\[ d\left(d_a\left\lf \ts\frac{y-1}{k-1} \right\rf+q_a y+1\right) - q_a\left(a\left\lf \ts\frac{y-1}{k-1} \right\rf+dy\right) = 
    \left\lf \ts\frac{y-1}{k-1} \right\rf + d > 0 \:\: \text{if}\:\: (k-1) \nmid y,  
\] 
and 
\[ (a+d)\left(d_a\left\lf \ts\frac{y-1}{k-1} \right\rf+q_a y+1\right) - (d_a+q_a)\left(a\left\lf \ts\frac{y-1}{k-1} \right\rf+dy\right) = 
    \left\lf \ts\frac{y-1}{k-1} \right\rf - y + a + d > 0 \:\: \text{if}\:\: (k-1) \mid y.
\]
Therefore,
\begin{equation} \label{lbasic3}
\max \left\{ f_{\ell}\left(\left\lf \ts\frac{y-1}{k-1} \right\rf,y\right): 1 \le y \le a-1 \right\} = f_{\ell}\left(\left\lf \ts\frac{a-2}{k-1} \right\rf,a-1\right). 
\end{equation}
The claim in eqn.~\eqref{lbasic1} follows from eqn.~\eqref{lbasic2} and eqn.~\eqref{lbasic3}. 
\end{Pf}

\begin{prop} \label{min}
Let $a,d,k$ be positive integers, with $\gcd(a,d)=1$ and $2 \le k \le a$. Let $S\big(AP(a,d;k)\big)=\la a,a+d,a+2d,\ldots,a+(k-1)d \ra$, and set $b=a+(k-1)d$. If $dd_b=1+q_b b$, with $1 \le d_b<b$ and $0 \le q_b<d$, and $a-2=\l(k-1)+\mu$, $0 \le \mu < k-1$, then 
\[ \min \left\{ \frac{g}{\left\lf \frac{d_b g}{b}\right\rf}: g \in \text{G}(S), d_b g \ge b \right\} = \frac{\l b+d}{\l d_b+q_b} =  
    \frac{\text{F}(S)-\mu d}{\left\lf \frac{d_b\big(\text{F}(S)-\mu d\big)}{b} \right\rf}, \]
with the minimum achieved at $g=\l a+d\big(\l(k-1)+1\big)=\l b+d$.
\end{prop}

\begin{Pf}
If $g \in \text{G}(S)$, then $g=ax+dy$, with $0 \le x \le \left\lf \frac{y-1}{k-1} \right\rf$ and $1 \le y \le a-1$. Let $d_b$ be such that $dd_b \equiv 1\pmod{b}$, $d_b>0$, and write $dd_b=q_b b+1$, $q_b \in {\Z}_{\ge 0}$. Then 
\begin{eqnarray*}
\frac{d_b g}{b} & = & \frac{d_b x\big(b-(k-1)d\big)+d_b dy}{b} \\
& = & d_b x + \frac{d_b d}{b}\big(y-(k-1)x\big) \\
& = & d_b x + q_b\big(y-(k-1)x\big) + \frac{y-(k-1)x}{b}. 
\end{eqnarray*} 
Since $0<y-(k-1)x<a<b$, we have  
\[ \frac{g}{\left\lf \frac{d_b g}{b}\right\rf} = \frac{ax+dy}{d_b x+q_b\big(y-(k-1)x\big)} \stackrel{\mbox{\tiny{def}}}{=} f_u(x,y). \]
We show that 
\begin{equation} \label{ubasic1}
\min \left\{ f_u(x,y): 0 \le x \le \left\lf \ts\frac{y-1}{k-1} \right\rf, 1 \le y \le a-1 \right\} = f_u\left( \l, \l(k-1)+1 \right)  
= \frac{\text{F}(S)-\mu d}{\left\lf \ts\frac{d^{-1}\big(\text{F}(S)-\mu d\big)}{b} \right\rf}. 
\end{equation}
\vskip 5pt

\noindent Fix $y \in \{1,\ldots,a-1\}$. For $x \in \left\{0,1,2,\ldots,\left\lf \ts\frac{y-1}{k-1} \right\rf -1\right\}$, the numerator of 
\[f_u(x+1,y) - f_u(x,y) = \frac{a(x+1)+dy}{d_b(x+1)+q_b\big(y-(k-1)(x+1)\big)} - \frac{ax+dy}{d_b x+q_b\big(y-(k-1)x\big)} \]
is 
\begin{eqnarray*} 
a\left(d_b x+q_b\big(y-(k-1)x\big)\right)-\big(d_b-q_b(k-1)\big)(ax+dy) = y\left( aq_b-d(d_b-q_b(k-1)\big) \right) \\
= y \left( q_b\big(a+(k-1)d\big) - (q_b b+1) \right) = -y < 0. 
\end{eqnarray*}
Therefore, 
\begin{equation} \label{ubasic2}
\min \left\{ f_u(x,y): 0 \le x \le \left\lf \ts\frac{y-1}{k-1} \right\rf, 1 \le y \le a-1 \right\} = \min \left\{ f_u\left(\left\lf \ts\frac{y-1}{k-1} \right\rf,y\right): 1 \le y \le a-1 \right\}. 
\end{equation}
\vskip 5pt

\noindent Since $y-(k-1)\left\lf \ts\frac{y-1}{k-1} \right\rf=1+\big(y-1 \!\!\pmod{k-1}\big)$, we may write   
\begin{equation} \label{ubasic3}
f_u\left(\left\lf \ts\frac{y-1}{k-1} \right\rf,y\right) = \frac{a\left\lf \ts\frac{y-1}{k-1} \right\rf+dy}{d_b\left\lf \ts\frac{y-1}{k-1} \right\rf + q_b\left(1+\big(y-1 \!\!\!\!\pmod{k-1}\big)\right)}. 
\end{equation}
\vskip 5pt

\noindent Write $y-1=q(k-1)+r$, $0 \le r<k-1$. Then the minimum on the right-side of eqn.~\eqref{ubasic2} can be written as
\begin{equation} \label{ubasic4}
\min \left\{ \frac{qa+dy}{qd_b+q_b(r+1)}: 0 \le q \le \l -1, 0 \le r < k-1 \;\text{or}\; q=\l, 0 \le r \le \mu \right\} 
\end{equation}
in view of eqn.~\eqref{ubasic3}. 
\vskip 5pt

\noindent Fix $q$. For $0 \le q<\l$, let $0 \le r < k-1$, and for $q=\l$, let $0 \le r \le \mu$. The numerator of 
\begin{eqnarray*}
f_u\left(\left\lf \ts\frac{y}{k-1} \right\rf,y+1\right) - f_u\left(\left\lf \ts\frac{y-1}{k-1} \right\rf,y\right) \\
= \frac{qa+d(y+1)}{qd_b+q_b(r+2)} - \frac{qa+dy}{qd_b+q_b(r+1)} = \frac{d\big(qd_b+q_b(r+1)\big) - q_b(qa+dy)}{\big(qd_b+q_b(r+1)\big) \big(qd_b+q_b(r+2) \big) } 
\end{eqnarray*}
is 
\[ q(q_b b-q_b a+1)-q_b d\big(y-(r+1)\big) = q + q\big(q_b(k-1)d\big) - q_b dq(k-1) = q > 0. \] 
\vskip 5pt

\noindent So for fixed $q$, the minimum is achieved when $r=0$, and we now have
\begin{eqnarray} \label{ubasic5}
\min \left\{ \frac{qa+dy}{qd_b+q_b(r+1)}: 0 \le q \le \l -1, 0 \le r < k-1 \;\text{or}\; q=\l, 0 \le r \le \mu \right\} \nonumber \\ 
= \min \left\{ \frac{qb+d}{qd_b+q_b}: 0 \le q \le \l \right\}.  
\end{eqnarray}
\vskip 5pt

\noindent Finally for $q \in \{0,1,2,\ldots,\l -1 \}$, the numerator of  
\[ \frac{(q+1)b+d}{(q+1)d_b+q_b} - \frac{qb+d}{qd_b+q_b} = \frac{b(qd_b+q_b) - d_b(qb+d)}{(qd_b+q_b) \big((q+1)d_b+q_b\big)} \]
is 
\[ q_b b - (q_b b+1) < 0. \]
\vskip 5pt

\noindent Therefore, the minimum is achieved at $q=\l$, and we have 
\begin{equation} \label{ubasic6}
\min \left\{ \frac{qb+d}{qd_b+q_b}: 0 \le q \le \l \right\} = \frac{\l b+d}{\l d_b+q_b} = f_u\left( \l, \l(k-1)+1 \right). 
\end{equation}
The claim in eqn.~\eqref{ubasic1} follows from eqns.~\eqref{ubasic2} - ~\eqref{ubasic6}. 
\end{Pf}

\begin{prop} \label{link}
Let $a,d$ be positive integers, with $\gcd(a,d)=1$ and let $b=a+(a-1)d$. Let $S\big(AP(a,d;a)\big)=\la a,a+d,a+2d,\ldots,a+(a-1)d\ra$, and let
\[ dd_a = 1+q_a a, \quad dd_b = 1+q_b b, \quad aa_b \equiv 1 \!\!\!\!\pmod{b}, \quad dd^{\pr} \equiv 1 \!\!\!\!\pmod{b-d}, \]
where 
\[ 1 \le d_a < a, \quad  0 \le q_a < d, \quad 1 \le d_b < b, \quad  0 \le q_b < d, \quad 1 \le a_b < b, \quad 1 \le d^{\pr} < b-d. \]
\begin{itemize}
\item[{\rm (i)}]
In terms of $d_a$ and $q_a$, we have 
\[ q_b = q_a, \quad d_b = d_a+(a-1)q_a, \quad d^{\pr} = d_a+(a-2)q_a, \quad a_b = (a-1)q_a+d_a+1. \]   
\item[{\rm (ii)}]
\[ \left\lc\frac{d_a g}{a}\right\rc = \left\lc\frac{a_b g}{b}\right\rc \;\;\text{and}\;\; \left\lf\frac{d_b g}{b}\right\rf = \left\lf \frac{d^{\pr} g}{b-d}
    \right\rf. 
\]
for each $g \in \text{G}(S)$. 
\end{itemize}
\end{prop}

\begin{Pf}
If $g \in \text{G}(S)$, then $g=dy$, with $1 \le y \le a-1$. From Propositions \ref{max} and \ref{min}, 
\begin{equation} \label{linkeq}
\left\lc\frac{d_a g}{a}\right\rc = q_a y + 1 \:\:\text{and}\:\: \left\lf\frac{d_b g}{b}\right\rf = q_b y. 
\end{equation}

Multiplying both sides of $dd_a-q_a a=1$ by $a-1$, then adding and subtracting $ad_a$, and rearranging, we get 
\[ a\big( (a-1)q_a + d_a + 1 \big) - bd_a = 1. \]
Thus, $a\big((a-1)q_a+d_a+1\big) \equiv 1\pmod{b}$ and $1 \le (a-1)q_a+d_a+1<b$, so that 
\begin{equation}\label{a_b}
a_b = (a-1)q_a + d_a  + 1. 
\end{equation}
Therefore
\[ \frac{a_b g}{b} = \frac{\big((a-1)q_a+d_a+1\big)dy}{b} = \frac{(b-a)q_a y+\big(1+q_a a\big)y+dy}{b} = q_a y + \frac{(d+1)y}{b}. \] 
Since $0 \le (d+1)y \le (d+1)(a-1) < b$, we have $\left\lc\frac{a_b g}{b}\right\rc=q_a y+1=\left\lc\frac{d_a g}{a}\right\rc$ by eqn.~\eqref{linkeq}. This proves the first part of the Proposition. 
\vskip 5pt

For $j \in \{0,\ldots,a-1\}$, we have $d\big(d_a+jq_a\big)-(a+jd)q_a=dd_a-q_a a=1$. Thus, $d\big(d_a+jq_a\big) \equiv 1\pmod{a+jd}$, and since $1 \le d_a+jq_a<a+jd$, we have $d_a+jq_a$ equals $d^{-1}\bmod{a+jd}$ for each $j \in \{0,\ldots,a-1\}$. In particular, 
\begin{equation} \label{d_b}
d_b = d_a+(a-1)q_a. 
\end{equation}
Multiplying both sides of eqn.~\eqref{d_b} by $d$ we get 
\[ 1+q_b b = dd_b = d\big(d_a+(a-1)q_a\big) = (1+q_a a) + (a-1)dq_a = 1+q_a b. \]  
Since $d^{\pr}$ is $d^{-1}\bmod{a+jd}$ for $j=a-2$, we have  
\begin{equation}\label{q_b}
d^{\pr} = d_a + (a-2)q_a \:\:\text{and}\:\: q_a = q_b.  
\end{equation}
Eqns.~\eqref{a_b}, \eqref{d_b}, \eqref{q_b} prove part (i). 
\vskip 5pt

We now have 
\[ \frac{d^{\pr} g}{b-d} = \frac{\big((a-2)q_a+d_a\big)dy}{b-d} = \frac{\big((b-d)-a\big)q_a y+\big(1+q_a a\big)y}{b-d} = q_b y + \frac{y}{b-d}. \] 
Since $0 \le y < b-d$, we have $\left\lf\frac{d^{\pr} g}{b-d}\right\rf=q_b y=\left\lf\frac{d_b g}{b}\right\rf$ by eqn.~\eqref{linkeq}. This proves part (ii). 
\end{Pf}
\vskip 5pt

\noindent Let $S=S\big(AP(a,d;k)\big)$ denote the numerical semigroup generated by $\text{AP}(a,d;k)$. Let $a,b,c$ be any positive integers. Since $S(a,b,c)=S(a\bmod{b},b,c)$ and $S(a,b,c)={\Z}_{\ge 0}$ for $c \ge a$, we may assume without loss of generality that $c<a<b$. Propositions \ref{basic3} and \ref{basic4} provide a one-one correspondence between $S(a,b,c)$ and $S([\a,\b])$. Via this correspondence, each interval $[\a,\b]$, with $\a \ge 0$, corresponds to a unique triple $(a,b,c)$, with $a \in \N$. The assumption $c<a<b$ on the positive integers $a,b,c$ gives rise to intervals $[\a,\b]$ with $\a>1$. In Theorem \ref{I}, we determine all intervals $[\a,\b]$ with $\a>1$ such that $S([\a,\b])=S$. We then use Proposition \ref{basic3} to determine all triples $(a,b,c)$ with $c<a<b$ such that $S(a,b,c)=S$. To extend the result of Theorem \ref{I} to all intervals $[\a,\b]$ with $\a \ge 0$, one may combine Proposition \ref{basic3} and the fact that $S(a,b,c)=S(a^{\pr},b,c)$ for $a \equiv a^{\pr}\pmod{b}$. 

We characterize intervals $[\a,\b]$, $\a>1$ in Theorem \ref{I} by calculating certain B\'{e}zout sequences related to $\text{AP}(a,d;k)$. Bullejos \& Rosales \cite{BR09} provide an algorithm to determine the B\'{e}zout sequence connecting any two rational numbers $p_1/q_1$ and $p_2/q_2$. However, it remains to determine $p_1/q_1$ and $p_2/q_2$ in order to characterize $[\a,\b]$. We give a direct method to explicitly determine the desired B\'{e}zout sequences.    
\vskip 5pt 

\begin{thm} \label{I}
Let $a,d,k$ be positive integers, with $\gcd(a,d)=1$, $2 \le k \le a$, $b=a+(k-1)d$, and let $\l=\left\lf\frac{a-2}{k-1}\right\rf$. Define $dd_a=1+q_a a$, 
where $1 \le d_a < a$ and $0 \le q_a < d$.
Let $S\big(AP(a,d;k)\big)=\la a,a+d,a+2d,\ldots,a+(k-1)d \ra$. Then for $\a,\b \in [1,\infty) \cap \Q$,
\[ S\big([\a,\b]\big) = S\big(AP(a,d;k)\big) \]
if and only if
\[ \frac{\l a+d(a-1)}{\l d_a+q_a(a-1)+1} < \a \le \frac{a}{d_a}, \quad \frac{b}{d_a+(a-1)q_a} \le \b < \frac{\l b+d}{\l\big(d_a+(a-1)q_a\big)+q_a}, \]
or
\begin{eqnarray*} 
\frac{\l b+d}{\l\big(b-d_a-(a-1)q_a\big)+(d-q_a)} < \a \le \frac{b}{b-d_a-(a-1)q_a},  \\
\frac{a}{a-d_a} \le \b < \frac{\l a+d(a-1)}{\l(a-d_a)+(d-q_a)(a-1)-1}. 
\end{eqnarray*}
\noindent Moreover, for $k=a$, we also have 
\[ \frac{d(a-1)}{q_a(a-1)+1} < \a \le \frac{b}{(a-1)q_a+d_a+1}, \quad \frac{b-d}{(a-2)q_a+d_a} \le \b < \frac{d}{q_a}, \]
or
\[ \frac{d}{d-q_a} < \a \le \frac{b-d}{b-d-(a-2)q_a-d_a}, \quad \frac{b}{b-(a-1)q_a-d_a-1} \le \b < \frac{d(a-1)}{(d-q_a)(a-1)-1}, \]
\noindent We adopt the convention $\frac{1}{0}=\infty$. 
\end{thm}

\begin{Pf}
Throughout this proof, we define
\[ dd_b = 1+q_b b, \quad aa_b \equiv 1 \!\!\!\!\pmod{b}, \quad dd^{\pr} \equiv 1 \!\!\!\!\pmod{b-d}, \]
with $1 \le d_b < b$, $0 \le q_b < d$, $1 \le a_b<b$, and $1 \le d^{\pr}<b-d$. 
\vskip 5pt

Let $\a,\b \in [1,\infty) \cap \Q$ be such that $S(I)=S(\text{AP}(a,d;k))=S$, with $I=[\a,\b]$. By Proposition \ref{basic1}, it is enough to show that $\frac{i}{n} \notin I$ when $i \in \text{G}(S)$ and $n \in \N$, and $\frac{j}{n} \in I$ when $j \in \text{AP}(a,d;k)$ for some $n=n(j) \in \N$ that depends on $j$. By eqn.~\eqref{Bez1}, there exists a proper B\'{e}zout sequence $\left\{\frac{a_1}{b_1},\ldots,\frac{a_k}{b_k}\right\}$ such that $\a \le \frac{a_1}{b_1}<\frac{a_k}{b_k} \le \b$ and $a_1,\ldots,a_k$ is a permutation of $\text{AP}(a,d;k)$. Then there exists $r \in \{1,\ldots,k\}$ such that 
\[ a_1 > \cdots > a_r < a_{r+1} < \cdots < a_k \]
by eqn.~\eqref{convex}. Therefore, $a_r=a$ and $a+d=\min\{a_{r-1},a_{r+1}\}$, unless $r=1$ or $k$. 

If $r=1$, then $a_1,\ldots,a_k$ is the increasing sequence $a,a+d,\ldots,a+(k-1)d$. Thus, $a_i=a+(i-1)d$, $1 \le i \le k$. The sequence of $b_i$'s satisfy $a_{i+1}b_i-a_ib_{i+1}=1$, so that $b_i \equiv a_{i+1}^{-1} \pmod{a_i} \equiv d^{-1} \pmod{a_i}$ for $1 \le i \le k-1$. From $a_kb_{k-1}-a_{k-1}b_k=1$ we have $b_k \equiv -(a_{k-1})^{-1} \pmod{a_k} \equiv d^{-1} \pmod{a_k}$. Thus, $b_i \equiv d^{-1} \pmod{a_i}$ for each $i \in \{1,\ldots,k\}$. Since $b_i,d^{-1}\pmod{a_i}$ belong to $\{1,\ldots,a_i-1\}$, we have $b_i=d^{-1}\pmod{a_i}$, for each $i \in \{1,\ldots,k\}$. So $\a,\b$ must satisfy  
\begin{equation} \label{r=1}
\a \le \frac{a_1}{b_1} = \frac{a}{d_a} < \frac{a_k}{b_k} = \frac{b}{d_b} \le \b. 
\end{equation}
In order that $\frac{i}{n} \notin I$ whenever $i \in \text{G}(S)$ and $n \in \N$, we must have either $\frac{i}{n}<\a$ or $\frac{i}{n}>\b$. From eqn.~\eqref{r=1}, this amounts to $\frac{i}{n}<\frac{a}{d_a}$ or $\frac{i}{n}>\frac{b}{d_b}$, or to $n>\frac{id_a}{a}$ or $n<\frac{id_b}{b}$. Therefore, $\a$ is greater than the maximum of $i/\left\lc\frac{id_a}{a}\right\rc$ as $i$ runs through values in $\text{G}(S)$ and $\b$ is less than the minimum of $i/\left\lf\frac{id_b}{b}\right\rf$ as $i$ runs through values in $\text{G}(S)$. Therefore, by Propositions \ref{max} and \ref{min} 
\begin{equation} \label{r=1bounds}
\a > \frac{\l a+d(a-1)}{\l d_a+q_a(a-1)+1} \;\;\text{and}\;\; \b < \frac{\l b+d}{\l d_b+q_b}. 
\end{equation}
This gives the first of the two ranges for $\a$ and $\b$ in the theorem using Proposition \ref{link}, part (i).  
\vskip 5pt

If $r=k$, then $a_1,\ldots,a_k$ is the decreasing sequence $a+(k-1)d,a+(k-2)d,\ldots,a$.  Thus, $a_i=a+(k-i)d$, $1 \le i \le k$. Arguing as above, we see that $b_i=a_i-\left(d^{-1} \pmod{a_i}\right)$ for each $i \in \{1,\ldots,k\}$. So $\a,\b$ must satisfy 
\begin{equation} \label{r=k}
\a \le \frac{a_1}{b_1} = \frac{b}{b-d_b} < \frac{a_k}{b_k} = \frac{a}{a-d_a} \le \b. 
\end{equation}
In order that $\frac{i}{n} \notin I$ whenever $i \in \text{G}(S)$ and $n \in \N$, we must have either $\frac{i}{n}<\a$ or $\frac{i}{n}>\b$. From eqn.~\eqref{r=k}, this amounts to $\frac{i}{n}<\frac{b}{b-d_b}$ or $\frac{i}{n}>\frac{a}{a-d_a}$, or to $n>\frac{i(b-d_b)}{b}$ or $n<\frac{i(a-d_a)}{a}$. Therefore, $\a$ is greater than the maximum of $i/\left\lc\frac{i(b-d_b)}{b}\right\rc$ as $i$ runs through values in $\text{G}(S)$ and $\b$ is less than the minimum of $i/\left\lf\frac{i(a-d_a)}{a}\right\rf$ as $i$ runs through values in $\text{G}(S)$. Since
\[ \frac{m}{\left\lc\frac{(b-d_b)m}{b}\right\rc} < \frac{n}{\left\lc\frac{(b-d_b)n}{b}\right\rc}  \LIff  
    \frac{n}{\left\lf\frac{d_b n}{b}\right\rf} < \frac{m}{\left\lf\frac{d_b m}{b}\right\rf} 
\]
and
\[ \frac{m}{\left\lf\frac{(a-d_a)m}{a}\right\rf} < \frac{n}{\left\lf\frac{(a-d_a)n}{a}\right\rf}  \LIff  
    \frac{n}{\left\lc\frac{d_a n}{a}\right\rc} < \frac{m}{\left\lc\frac{d_a m}{a}\right\rc} 
\]
for every $m,n \in \N$, we have 
\begin{equation} \label{r=kbounds}
\a > \frac{\l b+d}{\l(b-d_b)+(d-q_b)} \;\;\text{and}\;\; \b < \frac{\l a+d(a-1)}{\l(a-d_a)+(d-q_a)(a-1)-1}, 
\end{equation}
by Propositions \ref{max} and \ref{min}. This gives the second of the two ranges for $\a$ and $\b$ in the theorem using Proposition \ref{link}, part (i).
\vskip 5pt

Henceforth, assume $r \ne 1,k$. Two cases arise: (i) $a_{r+1}=a+d$, and (ii) $a_{r-1}=a+d$. These cases arise only if $k=a$, as we show below. 
\vskip 5pt

\noindent {\sc Case} (i) From $a_r \mid (a_{r-1}+a_{r+1})$ we have $a_{r-1} \equiv -a_{r+1} \equiv -d\pmod{a}$. If $a_{r-1}=a+id$, with $0 \le i \le k-1 \le a-1$, then $i \equiv -1\pmod{a}$. This is possible only if $k=a$, and then $a_{r-1}=a+(a-1)d=b$. From eqn.~\eqref{convex}, the sequence $a_1,\ldots,a_k$ must be $b=a+(a-1)d,a,a+d,\ldots,a+(a-2)d$. The corresponding sequence $b_2,\ldots,b_k$ is the same as in the case $r=1$ above. Moreover, $ab_1-\big(a+(a-1)d\big)b_2=1$, so that $b_1 \equiv a_b \pmod{b}$. Since $b_1,a_b$ belong to $\{1,\ldots,b-1\}$, we have $b_1=a_b$. 
So $\a,\b$ must satisfy 
\begin{equation} \label{r=1,k}
\a \le \frac{a_1}{b_1} = \frac{b}{a_b} < \frac{a_k}{b_k} = \frac{b-d}{d^{\pr}} \le \b. 
\end{equation}
In order that $\frac{i}{n} \notin I$ whenever $i \in \text{G}(S)$ and $n \in \N$, we must have either $\frac{i}{n}<\a$ or $\frac{i}{n}>\b$. From eqn.~\eqref{r=1,k}, this amounts to $\frac{i}{n}<\frac{b}{a_b}$ or $\frac{i}{n}>\frac{b-d}{d^{\pr}}$, or to $n>\frac{i a_b}{b}$ or $n<\frac{i d^{\pr}}{b-d}$. Therefore, $\a$ is greater than the maximum of $i/\left\lc\frac{i a_b}{b}\right\rc$ as $i$ runs through values in $\text{G}(S)$ and $\b$ is less than the minimum of $i/\left\lf\frac{i d^{\pr}}{b-d}\right\rf$ as $i$ runs through values in $\text{G}(S)$. By Proposition \ref{link}, part (ii), the bounds for $\a$ and $\b$ are as given by eqn.~\eqref{r=1bounds} with $\l=0$. This gives the first of the two ranges for $\a$ and $\b$ in the additional case $k=a$ in the theorem using Proposition \ref{link}, part (i). 
\vskip 5pt

\noindent {\sc Case} (ii) From $a_r \mid (a_{r-1}+a_{r+1})$ we have $a_{r+1} \equiv -a_{r-1} \equiv -d\pmod{a}$. Arguing as in Case (i), the sequence $a_1,\ldots,a_k$ must be $a+(a-2)d,a+(a-3)d,\ldots,a,a+(a-1)d=b$, and this is possible only if $k=a$. The corresponding sequence $b_1,\ldots,b_{k-1}$ is the same as in the case $r=k$ above, and $b_k=b-\left(a_b \pmod{b}\right)$. So $\a,\b$ must satisfy 
\begin{equation} \label{r=k}
\a \le \frac{a_1}{b_1} = \frac{b-d}{b-d-d^{\pr}} < \frac{a_k}{b_k} = \frac{b}{b-a_b} \le \b. 
\end{equation}
In order that $\frac{i}{n} \notin I$ whenever $i \in \text{G}(S)$ and $n \in \N$, we must have either $\frac{i}{n}<\a$ or $\frac{i}{n}>\b$. From eqn.~\eqref{r=k}, this amounts to $\frac{i}{n}<\frac{b-d}{b-d-d^{\pr}}$ or $\frac{i}{n}>\frac{b}{b-a_b}$, or to $n>\frac{i(b-d-d^{\pr})}{b-d}$ or $n<\frac{i(b-a_b)}{b}$. Therefore, $\a$ is greater than the maximum of $i/\left\lc\frac{i(b-d-d^{\pr})}{b-d}\right\rc$ as $i$ runs through values in $\text{G}(S)$ and $\b$ is less than the minimum of $i/\left\lf\frac{i(b-a_b)}{b}\right\rf$ as $i$ runs through values in $\text{G}(S)$. Since
\[ \frac{m}{\left\lc\frac{(b-d-d^{\pr})m}{b-d}\right\rc} < \frac{n}{\left\lc\frac{(b-d-d^{\pr})n}{b-d}\right\rc}  \LIff  
     \frac{n}{\left\lf\frac{d^{\pr}n}{b-d}\right\rf} < \frac{m}{\left\lf\frac{d^{\pr}m}{b-d}\right\rf} 
\]
and
\[ \frac{m}{\left\lf\frac{(b-a_b)m}{b}\right\rf} < \frac{n}{\left\lf\frac{(b-a_b)n}{b}\right\rf}  \LIff  
     \frac{n}{\left\lc\frac{a_b n}{b}\right\rc} < \frac{m}{\left\lc\frac{a_b m}{b}\right\rc} 
\]
for every $m,n \in \N$, we have the bounds for $\a$ and $\b$ as given by eqn.~\eqref{r=kbounds} with $\l=0$ by Proposition \ref{link}, part (ii). This gives the second of the two ranges for $\a$ and $\b$ in the additional case $k=a$ in the theorem using Proposition \ref{link}, part (i).
\end{Pf}
\vskip 5pt

\begin{cor} \label{spcase1}
Let $a,k$ be positive integers, with $2 \le k \le a$, $b=a+k-1$, and $\l=\left\lf\frac{a-2}{k-1}\right\rf$. Let $S\big(AP(a,1;k)\big) = \la a,a+1,a+2,\ldots,a+k-1 \ra$. Then for $\a,\b \in [1,\infty) \cap \Q$,
\[ S\big([\a,\b]\big) = S\big(AP(a,1;k)\big) \]
if and only if
\[ a-\frac{1}{\l+1} < \a \le a, \quad b \le \b < b+\frac{1}{\l} \]
or
\[ 1+\frac{\l}{\l(b-1)+1} < \a \le 1+\frac{1}{b-1}, \quad 1+\frac{1}{a-1} \le \b < 1+\frac{\l+1}{(\l+1)(a-1)-1}. \]
\noindent Moreover, for $k=a$, we also have 
\[ a-1 < \a \le \frac{b}{2}, \quad b-1 \le \b, \]
or
\[ 1 < \a \le 1+\frac{1}{b-2}, \quad 1+\frac{2}{b-2} \le \b < 1+\frac{1}{a-2}. \] 
\noindent We adopt the convention $\frac{1}{0}=\infty$. 
\end{cor}
\vskip 5pt

\begin{rem} \label{spcase1_k=a}
For any positive integer $a>1$ and for $\a,\b \in [1,\infty) \cap \Q$, 
\[ S\big([\a,\b]\big) = S\big(AP(a,1;a)\big) = \big\{0,a,\to\big\} \]
is the special case $\l=0$ in Corollary \ref{spcase1}. We adopt the convention $\frac{1}{0}=\infty$. 
\end{rem}
\vskip 5pt

\begin{thm} \label{II}
Let $a,d,k$ be positive integers, with $\gcd(a,d)=1$, $2 \le k \le a$, $b=a+(k-1)d$, and let $\l=\left\lf\frac{a-2}{k-1}\right\rf$. Define $dd_a=1+q_a a$, 
where $1 \le d_a < a$ and $0 \le q_a < d$. Let $S\big(AP(a,d;k)\big)=\la a,a+d,a+2d,\ldots,a+(k-1)d \ra$. Then for positive integers $p,q,m$ with 
$q<p<m$, 
\[ S\big(p,m,q\big) = S\big(AP(a,d;k)\big) \]
if and only if
\begin{eqnarray*}
\frac{md_a}{a} \le p < \frac{m\big(\l d_a+q_a(a-1)+1\big)}{\l a+d(a-1)}, \\
p + \frac{m\big((a-1)q_a+d_a\big)}{b} \le q < p - \frac{m\left(\l\big((a-1)q_a+d_a\big)+q_a\right)}{\l b+d}
\end{eqnarray*}
or
\begin{eqnarray*} 
m - \frac{m\big((a-1)q_a+d_a\big)}{b} \le p < m - \frac{m\left(\l\big((a-1)q_a+d_a\big)+q_a\right)}{\l b+d}, \\ 
p - m + \frac{md_a}{a} \le q < p - m - \frac{m\big((\l d_a+q_a)(a-1)-1\big)}{\l a+d(a-1)}.
\end{eqnarray*}
\noindent Moreover, for $k=a$, we also have 
\[ \frac{m\big((a-1)q_a+d_a+1\big)}{b} \le p < \frac{m\big(q_a(a-1)+1\big)}{d(a-1)}, \quad p - \frac{m\big((a-2)q_a+d_a\big)}{b-d} \le q < p - \frac{mq_a}{d}, \]
or
\begin{eqnarray*} 
m - \frac{m\big((a-2)q_a+d_a\big)}{b-d} \le p < m - \frac{mq_a}{d}, \\ 
p - m + \frac{m\big((a-1)q_a+d_a+1\big)}{b} \le q < p - m + \frac{m\big(q_a(a-1)+1\big)}{d(a-1)}.  
\end{eqnarray*}
\end{thm}

\begin{Pf}
This follows directly from Proposition \ref{basic2} and Theorem \ref{I}. If $S\big(p,m,q\big)=S\big([\a,\b]\big)$, then $\frac{m}{p}=\a$ and                $\frac{m}{p-q}=\b$, so that $p=\frac{m}{\a}$ and $p-q=\frac{m}{\b}$. The latter implies $q=p-\frac{m}{\b}$. The bounds on $\a,\b$ in Theorem \ref{I} translate to bounds on $p,q$ in terms of $m$ in Theorem \ref{II}.  
\end{Pf}
\vskip 5pt

\begin{cor} \label{spcase2}
Let $a,k$ be positive integers, with $2 \le k \le a$, $b=a+k-1$, and $\l=\left\lf\frac{a-2}{k-1}\right\rf$. Let $S\big(AP(a,1;k)\big) = \la a,a+1,a+2,\ldots,a+k-1 \ra$.  Then for positive integers $p,q,m$ with $q<p<m$, 
\[ S\big(p,m,q\big) = S\big(AP(a,1;k)\big) \]
if and only if
\[ \frac{m}{a} \le p < \frac{m\big(\l+1\big)}{(\l+1)a-1}, \quad p + \frac{m}{b} \le q < p - \frac{m\l}{\l b+1} \]
or
\[ m - \frac{m}{b} \le p < m - \frac{m\l}{\l b+1}, \quad p - m + \frac{m}{a} \le q < p - m - \frac{m\big(\l(a-1)-1\big)}{(\l+1)a-1)}. \]
\noindent Moreover, for $k=a$, we also have 
\[ \frac{2m}{b} \le p < \frac{m}{a-1}, \quad p - \frac{m}{b-1} \le q < p, \]
or
\[ m - \frac{m}{b-1} \le p < m, \quad p - m + \frac{2m}{b} \le q < p - m + \frac{m}{a-1}. \]
\end{cor}
\vskip 5pt

\begin{rem} \label{spcase2_k=a}
For any positive integer $a>1$ and positive integers $p,q,m$ with $q<p<m$,  
\[ S\big(p,m,q\big) = S\big(AP(a,1;a)\big) = \big\{0,a,\to\big\} \]
is the special case $\l=0$ in Corollary \ref{spcase2}. 
\end{rem}


\vskip 20pt

\begin{thebibliography}{99}

\bibitem{Bat58}
P. T. Bateman, Remark on a recent note on linear forms, {\it Amer. Math. Monthly\/} {\bf 65} (1958), 517--518. 

\bibitem{Bra42}
A. Brauer, On a problem of partitions, {\it Amer. J. Math.\/} {\bf 64} (1942), 299--312. 

\bibitem{BR09}
M. Bullejos and J. C. Rosales, Proportionally modular Diophantine inequalities and the Stern-Brocot tree, {\it Math. Comp.\/} {\bf 78} (2009), 1211--1226. 

\bibitem{DR06}
M. Delgado and J. C. Rosales, On the Frobenius number of a proportionally modular Diophantine inequality, {\it Port. Math.\/} {\bf 63} (2006), 
415--425. 

\bibitem{Gra73}
D. D. Grant, On linear forms whose coeffficients are in arithmetic progression, {\it Israel J. Math.\/} {\bf 15} (1973), 204--209.

\bibitem{NW72}
A. Nijenhuis and H. S. Wilf, Representations of integers by linear forms in nonnegative integers, {\it J. Number Theory\/} {\bf 4} (1972), 98--106. 

\bibitem{Rob56}
J. B. Roberts, Note on linear forms, {\it Proc. Amer. Math. Soc.\/} {\bf 7} (1956), 465--469. 

\bibitem{RGGU03}
J. C. Rosales, P. A. Garc\'{i}a-S\'{a}nchez, J. I. Garc\'{i}a-Garc\'{i}a and J. M. Urbano-Blanco, Proportionally Modular Diophantine Inequalities, {\it J. Number Theory\/} {\bf 103} (2003), 281--294. 

\bibitem{RG-S08}
J. C. Rosales and P. A. Garc\'{i}a-S\'{a}nchez, Numerical semigroups having a Toms decomposition, {\it Canad. Math. Bull.\/} {\bf 51} (1) (2008), 134--139.  

\bibitem{RGU05}
J. C. Rosales, P. A. Garc\'{i}a-S\'{a}nchez and J. M. Urbano-Blanco, Modular Diophantine inequalities and numerical semigroups, {\it Pacific J. Math.\/} {\bf 218} (2005), 379--398. 

\bibitem{RGU08}
J. C. Rosales, P. A. Garc\'{i}a-S\'{a}nchez and J. M. Urbano-Blanco, The set of solutions of a proportionally modular Diophantine inequality, {\it J. Number Theory\/} {\bf 128} (2008), 453--467. 


\bibitem{Tom03}
A. Toms, Strongly perforated $K_0$-groups of simple $C^{\star}$-algebras, {\it Canad. Math. Bull.\/} {\bf 46} (3) (2003), 457--472.  

\bibitem{Tri94}
A. Tripathi, The coin exchange problem for arithmetic progressions, {\it Amer. Math. Monthly\/} {\bf 101} No.10 (1994), 779--781.

\end{thebibliography}
\end{document}